\documentclass[11pt,3p]{scrartcl}
 
\usepackage[hidelinks]{hyperref}
\usepackage[pdftex,x11names]{xcolor}
\usepackage{subfigure}
\usepackage[T1]{fontenc}
\usepackage[utf8]{inputenc}
\usepackage{lmodern}
\usepackage[english]{babel}
\usepackage{amsmath}
\usepackage{amssymb}
\usepackage{amsfonts}
\usepackage{amsthm}
\usepackage[centerdot]{mathtools}
\usepackage{dsfont}
\newcommand{\R}{\mathds{R}}
\newcommand{\Z}{\mathds{Z}}

\usepackage{pifont}%
\newcommand{\cmark}{\text{\ding{51}}}%
\newcommand{\xmark}{\text{\ding{55}}}%

\DeclareMathOperator{\conv}{conv}

\usepackage{booktabs}

\usepackage{hyperref}
\usepackage[ruled,noend]{algorithm2e}
\SetKw{Continue}{continue}
\SetKw{Null}{null}
\SetKw{Break}{break}

\SetCommentSty{mycommfont}

\usepackage[textsize=tiny]{todonotes}
\usepackage[normalem]{ulem}
\usepackage{float}
\usepackage{authblk}

\usepackage{tikz}
\usepackage{tikz,tikz-3dplot}

\tikzset{
    cross/.pic = {
    \draw[rotate = 45] (-#1,0) -- (#1,0);
    \draw[rotate = 45] (0,-#1) -- (0, #1);
    }
}
\usetikzlibrary{shapes.misc}
\usetikzlibrary{patterns}
\usetikzlibrary{graphs,quotes}
 \usepackage{pgfplots}
 \pgfplotsset{compat=1.16,
    every axis/.append style={
        axis lines=center,
        xlabel style={anchor=south west},
        ylabel style={anchor=south west},
        zlabel style={anchor=south west},
        tick align=outside,}
}
 \usepgfplotslibrary{patchplots}
\tdplotsetmaincoords{60}{120}

\usepackage{xpatch}
\newcommand{\strictly}{}
\newcommand{\weakly}{weakly }
\newcommand{\X}{\mathcal{X}}

\newcommand{\Y}{\mathcal{Y}}

\usepackage[]{cleveref} %nameinlink makes all blue

\newtheoremstyle{thm}% name
{15pt}% Space above
{5pt}% Space below 
{\itshape}% Body font
{}% Indent amount: Indent amount: empty = no indent, \parindent = normal paragraph indent
{\bfseries}% Theorem head font\part{title}
{}% Punctuation after theorem head
{0.5em}% Space after theorem head: Space after theorem head: { } = normal interword space; \newline = linebreak
{}% Theorem head spec (can be left empty, meaning `normal')
\theoremstyle{thm}% default

\usepackage{xpatch}

\xpatchcmd{\proof}{\itshape}{\prooflabelfont}{}{}
\newcommand{\prooflabelfont}{\bfseries}

\newtheorem{thm}{Theorem}[section]
\newtheorem{lem}[thm]{Lemma}
\newtheorem{definition}[thm]{Definition}

\newtheorem{example}[thm]{Example}

\title{An Output-Polynomial Time Algorithm to Determine all Supported Efficient Solutions for Multi-Objective Integer Network Flow Problems}

\author[1]{David K\"onen}
\author[1]{Michael Stiglmayr}
\affil[1]{%
	University of Wuppertal\\
	School of Mathematics and Natural Sciences\\
	Optimization Group\\
	Gaußstraße 20, 42103 Wuppertal, Germany\\ 
}
\affil[]{E-Mail:~\href{mailto:koenen@uni-wuppertal.de}{koenen@uni-wuppertal.de}, \href{mailto:stiglmayr@uni-wuppertal.de}{stiglmayr@uni-wuppertal.de} }

\date{}

\begin{document}

\maketitle

	\begin{abstract}\small
	%% Text of abstract
	This paper addresses the problem of enumerating all \emph{supported efficient solutions} for a linear multi-objective integer minimum cost flow problem (MOIMCF). 
   It derives an output-polynomial time algorithm to determine all supported efficient solutions for MOIMCF problems. 
    This is the first approach to solve this general problem in output-polynomial time. 
    Moreover, we prove that the existence of an output-polynomial time algorithm to determine all weakly supported nondominated vectors (or all weakly supported efficient solutions) for a MOIMCF problem with a fixed number of $d\geq 3$ objectives can be excluded unless  $\mathbf{P} = \mathbf{N P}$. 
	
\end{abstract}

%%Graphical abstract
%\begin{graphicalabstract}
%\includegraphics{grabs}
%\end{graphicalabstract}

%%Research highlights
%\begin{highlights}
%\item Research highlight 1
%\item Research highlight 2
%\end{highlights}
	
	\par\vskip\baselineskip\noindent

\section{Introduction}

The minimum cost integer network flow problem is a fundamental and well-studied problem in combinatorial optimization \cite{Ahuja1993,bertsekas98network}.  The problem involves moving flow through a network at minimal cost. The given network has specific nodes that supply the flow units while others demand it. The flow units are moved along the arcs of the network. The flow along these arcs is constrained by specified lower and upper bounds, known as capacities, and has associated costs per flow unit. The objective is to find a flow that satisfies supply, demand, and capacity constraints and minimizes the overall cost. The minimum cost flow problem has numerous real-world applications, especially in industry and decision-making, such as inventory planning, data scaling, lot sizing, location problems, DNA sequence alignment, and project management~\cite{Ahuja1993}. Moreover, the minimum cost flow problem contains several combinatorial optimization problems as special cases. Among others, the important \emph{transportation problem}, and the \emph{assignment problem} reduce to the minimum cost flow problem.

For the single-objective version of this problem, there are various polynomial algorithms. A comprehensive overview of the related literature is provided in~\cite{Ahuja1993}. However, real-world problems often involve multiple conflicting objectives and no solution optimizes all objectives simultaneously.
These multi-objective problems with conflicting objectives arise in many real-world problems, including logistics, economics, finance, and more.  
In these conflicting scenarios, there is no solution that optimizes all objectives simultaneously. Thus, one is interested in finding solutions with the property that it can only be improved with respect to one objective if at least one other objective deteriorates. Such a solution is called \emph{efficient solution} or \emph{Pareto-optimal} solution and its image is called \emph{nondominated vector}.

Various approaches exist for tackling multi-objective optimization problems. The primary approaches can be categorized as follows. In \emph{a-priori methods}, the decision-maker specifies his/her preferences, by choosing weights or goals, before the optimization process begins. In the \emph{a-posteriori} methods, the optimization process generates a set of efficient solutions, allowing the decision maker to select the most suitable one. In \emph{interactive methods}, the decision-maker iteratively controls the optimization process, refining or adapting preferences and guiding the search to the most suitable solution.  

This paper focuses on a-posteriori methods, aiming to determine (all) or a suitable subset of the efficient solutions or nondominated vectors; for a summary of solution concepts in multi-objective optimization, see~\cite{serafini87some}.
One subset of interest is the set of all \emph{supported efficient solutions}, which are efficient solutions that can be obtained as optimal solutions of a \emph{weighted sum scalarization} with weights strictly greater than zero. 
\emph{Extreme supported solutions} are efficient solutions whose images lie on vertices of the upper image of the feasible outcome vectors, defined as the polyhedron given by the convex hull of feasible outcome vectors plus the (closed) positive orthant. While multi-objective linear optimization problems contain only supported efficient solutions, the efficient solution set of multi-objective combinatorial optimization problems, such as the multi-objective minimum cost integer network flow problem, also contains unsupported efficient solutions. The unsupported efficient solutions typically outnumber the supported ones, where the latter are easier to determine,  can serve as high-quality representations~\cite{Serpil2024}, and can be used as a foundation for two-phase methods to generate the entire nondominated point set. Despite their importance, several characterizations for supported efficient solutions and thus supported nondominated points are used in the literature as recently shown in~\cite{koenen2025supportednessmultiobjectivecombinatorialoptimization}. The paper follows the strict separation between supported and weakly supported solutions as presented therein.

\emph{Multi-objective integer minimum cost flow} (MOIMCF) problems are computationally much harder to solve than their single-objective counterparts. MOIMCF had been reviewed in~\cite{Hamacher2007}, where the authors comment on the lack of computationally efficient algorithms. However, there is still no efficient algorithm for determining all supported efficient solutions for MOIMCF problems. We give an overview of the new results below. 

  While its linear relaxation, the \emph{multi-objective minimum cost flow problem} (MOMCF), has only \emph{supported efficient solutions}, the \emph{multi-objective integer minimum cost flow} (MOIMCF) problem, on the other hand, can have supported, weakly supported, and unsupported efficient solutions.

There are several algorithms to determine efficient solutions for bi-objective integer minimum cost flow (BOIMCF) problems, e.g.,~\cite{Eusebio09,eusebio14,raith09,sedeno-noda01,sedona00,sedona03}. Raith and Sede\~no-Noda introduced an enhanced parametric approach to determine all extreme supported efficient solutions for BOIMCF problems~\cite{raith17}.
However, there are only quite a few specific methods designed to determine all (or subsets) of the nondominated vectors in the objective space (nor the corresponding efficient solutions in the decision space) for MOIMCF problems~\cite{EUSEBIO200968,fonesca10,sun11}.

Eus\'{e}bio and Figueria~\cite{EUSEBIO200968} introduced an algorithm that enumerates all supported nondominated vectors/efficient solutions for MOIMCF problems (assuming extreme supported solutions and corresponding weight vectors are given), based on zero-cost cycles in the incremental graph associated with the corresponding parametric network flow problems. They conclude that their proposed algorithm is the first step in developing further zero-cost cycle algorithms for solving MOIMCF problems.  However, they do not provide a specific method for determining those zero-cost cycles.  As mentioned before, different definitions exist for supported efficient solutions, which lead to different sets of supported efficient solutions, namely supported- and \strictly \weakly supported efficient solutions. Using the definition in~\cite{EUSEBIO200968} would also include the weakly supported efficient solutions. However, it can be shown that their proposed algorithm can only determine all supported efficient solutions. We refer to the example in~\Cref{2}. 

This paper will show that there is no \emph{output-polynomial} time algorithm to determine  all \weakly supported nondominated vectors (or \weakly supported efficient solutions) for a MOIMCF problem with a fixed number of objectives, unless $\mathbf{P} = \mathbf{N P}$. In contrast, we derive \emph{output-polynomial enumeration} algorithms to determine all \strictly supported efficient solutions, first for the bi-objective case and afterward for every fixed number of objectives. Here, output-polynomial time refers to a computation time that can be bounded by a polynomial in the input and output size.  We refer to the paper of Johnson, Papadimitriou and Yannakakis~\cite{JOHNSON1988119} for a detailed survey of output-sensitive complexity. 

The approach consists of two phases:  First, the algorithm determines all extreme supported points of the upper image and the associated weight vector for each maximally nondominated face. Hence, the approach successively determines all efficient solutions for each maximally nondominated face by determining all optimal solutions for the linear weighted-sum scalarization (single-objective parametric network flow) problem with the corresponding weight vector using the previously presented algorithm to determine all optimum integer flows in a network~\cite{KONEN2022333}. The method successively searches for \emph{proper zero-cost cycles} in linear time by using a modified depth-first search technique. 

 Given a BOIMCF problem and using the enhanced parametric network approach~\cite{raith17} to determine all $N$ extreme nondominated points in $\mathcal{O}(N\cdot n(m+n\log n ))$ time in a first step, results in an $\mathcal{O}(N\cdot n(m+n\log n ) + \mathcal{S}(m+n))$ time algorithm to determine all $\mathcal{S}$  supported efficient solutions to a BOIMCF problem.

Given a MOIMCF problem with a fixed number of $d$ objectives,  the dual Benson's algorithm~\cite{Ehrgott2012} can be used for the first phase. In addition, while the lexicographic MCFP can be solved in polynomial time~\cite{Hamacher2007}, we can use the lexicographic dual Benson's algorithm presented in \cite{Boekler2015}, which determines all extreme nondominated vectors, all facets of the \emph{upper image} and the \emph{weight space decomposition} in $\mathcal{O}(N^{\lfloor \frac{d}{2}\rfloor}(\operatorname{poly}(n,m)+ N \log N))$. Note that in phase two, the corresponding image of some solutions may lie in more than one maximally nondominated face and, therefore, are included in more than one face. The approach yields an $\mathcal{O}(N^{\lfloor \frac{d}{2}\rfloor}( \operatorname{poly}(n,m) + N\log N + N^{\lfloor \frac{d}{2}\rfloor}\mathcal{S}( m+n + N^{\lfloor \frac{d}{2}\rfloor}d))$ time algorithm to determine all $\mathcal{S}$ \strictly supported efficient solutions for a $d$-objective MOIMCF problem. To our knowledge, this is the first output-polynomial time algorithm to determine the complete set of supported efficient solutions.

The following table summarizes the contribution of our article and the existing results from the literature on the existence of output-polynomial time algorithms for the MOIMCF problem w.r.t.\ the different subsets of interest. 
A check mark indicates existence, a cross indicates that the existence of such an algorithm can be ruled out unless $\mathbf{P}= \mathbf{NP}$, and the question mark indicates that this problem remains an open question.

\begin{table}[htb]\footnotesize
	\begin{center}
        \caption{Results on the existence of output-polynomial time algorithms for the MOIMCF problem w.r.t.\ the different subsets of interest.}\label{tab1}
		\begin{tabular}{p{7em}|cccc}
			\toprule
			& \parbox{5em}{extreme supported}\rule[-1em]{0pt}{2.5em} & supported & \parbox{5em}{weakly supported} & \text {all} \\\midrule
			nondominated\linebreak vectors & \cmark\footnotemark[1] &   \textbf{?} & \xmark\footnotemark[2]  &  \xmark\footnotemark[3] \\\midrule
			efficient\linebreak solutions & \cmark\footnotemark[4] & \cmark\footnotemark[5] & \xmark\footnotemark[2]  & \xmark\footnotemark[6] \\
			\bottomrule
			\multicolumn{4}{p{8cm}}{\scriptsize
				\textsuperscript{1}\; \cite{Boekler2015,Ehrgott2012}\rule{0pt}{4ex} \hspace{1.5cm} \textsuperscript{4}\; Transferable from \cite{Boekler2015,Ehrgott2012} \newline
				\textsuperscript{2}\; \Cref{thm:moimcf} \hspace{0.51cm} \textsuperscript{5}\; \Cref{thm:outputpoly} \newline
				\textsuperscript{3}\;  \cite{Boekler2017} \hspace{1.83cm} \textsuperscript{6}\; Transferable from \cite{Boekler2017}, see~\Cref{lem:mosp_x}
			}
		\end{tabular}
	\end{center}
\end{table}

Note that an output-polynomial time algorithm for determining all supported efficient solutions is insufficient for the computation of all nondominated supported vectors in output-polynomial time since there may be exponentially many solutions mapping to the same vectors. It is easy to check that all extreme supported efficient solutions can be determined in output-polynomial time due to the existence of the output-polynomial time algorithm for the extreme supported nondominated vectors. Note that the algorithm of Eus\'{e}bio and Figueria~\cite{EUSEBIO200968} would also determine all supported efficient solutions in output-polynomial time when all maximally nondominated faces and corresponding weight vectors are already given.

The remainder of the paper is structured as follows. \Cref{2} presents the notations and preliminaries, proves that there exist weakly supported solutions in multi-objective minimum integer cost flow problems, and proves that there does not exist an output-polynomial time algorithm to determine all minimum cost flow problems in output-polynomial time unless  $\mathbf{P} = \mathbf{N P}$. \Cref{3} presents the algorithms to determine all supported efficient solutions for BOIMCF problems as well as all \strictly supported efficient solutions for MOIMCF problems. \Cref{chapt:concl} summarizes the article and gives an outlook on future research.

\section{Multi-Objective Integer Network Flows}
\label{2}
This section formally introduces the multi-objective integer minimum cost flow problem along with some of the most important results and properties of network flows. 
For a comprehensive introduction to the graph-theoretic foundations and the basics of network flow theory, we refer to \cite{Ahuja1993,Diestel2006}.

Let  $\R^d_{\geqq} \coloneqq  \{x\in \R^d \colon x \geqq 0\} $ denote the non-negative orthant of $\R^d$. $\R^d_{\geq}$ and its interior $\R^d_{>}$ is defined accordingly. In our notation, we also use the Minkowski sum and the Minkowski product of two sets $A, B \subseteq \R^{d}$, which are defined as $A+B\coloneqq \{a+b \colon a \in A, b \in B\}$ and $A \cdot B\coloneqq \{a \cdot b \colon $ $a \in A, b \in B\}$, respectively.

Furthermore, let $D=(V,A)$ be a directed graph with $|V|=n$ nodes and $|A| = m$ arcs 
together with integer-valued, non-negative, finite lower and upper capacity bounds $l_{ij}$ and $u_{ij}$, respectively, for each arc $(i,j)\in A$. 
Moreover, let $b_i$ be the integer-valued \emph{flow balance} of the node $i\in V$ and let $c\colon\, \R^m \rightarrow \R^d_{\geqq}$ be the (vector valued) objective function and $C\in\R^{d\times m}$ the corresponding cost matrix, i.\,e., $c(f)=C\cdot f$. 
 Then the \emph{multi-objective integer minimum cost flow} (MOIMCF) problem can be formulated as:
\[\def\arraystretch{1.4}
  \begin{array}{rr@{\extracolsep{0.7ex}}c@{\extracolsep{0.7ex}}lll}
    \min  &\multicolumn{1}{l}{ c(f) }\\
    \text{s.t.}  &\displaystyle \sum_{j:(i,j)\in A}f_{ij} - \sum_{j:(j,i)\in A} f_{ji} &=& b_i 	\quad  & i\in V, & \text{
    (\emph{flow balance constraint})}\\
     &l_{ij} \leq f_{ij} &\leq& u_{ij} \quad   & (i,j)\in A, &\text{
    (\emph{capacity constraint}) }\\ 
    & f_{ij} &\in& \Z_{\geq}  & (i,j)\in A.
   \end{array}
\]

Where we denote a feasible solution $f$ as \emph{flow} and define the \emph{cost} of a flow $f$ as  
\[
  c(f)\coloneqq \bigl(c^1(f),\dots,c^d(f)\bigr)^\top  \coloneqq  \Bigl(\sum_{(i,j) \in A} c^1_{ij}\, f_{ij}, \dots, \sum_{(i,j) \in A} c^d_{ij} \,f_{ij}\Bigr)^\top =C\cdot f.
\]

The \emph{(continuous) multi-objective minimum cost flow} (MOMCF) problem can be formulated as the LP-relaxation of MOIMCF. From now on, we assume that $D$ is connected  (see~\cite{Ahuja1993}) and that there is at least one feasible flow.

Then, we denote the set of feasible outcome vectors in objective space by $\mathcal{Y} \coloneqq   \{C\, f \colon f\in \mathcal{X} \}$, where $\mathcal{X}$ is the set of all feasible flows.
We assume that the objective functions are conflicting, which implies that we exclude the existence of an \emph{ideal} solution that minimizes all objectives simultaneously. 
Throughout this article, we will use the \emph{Pareto concept of optimality}, which is based on the componentwise order: 
\[
c(f)  \leqslant c(f^{\prime}) \iff c^k(f) \leqq c^k(f^{\prime}) \quad \forall k \in \{1,\dots , d\} \quad \text { and } \quad c(f) \neq c(f^{\prime})
\]

A vector $c(f') \in \R^{d}$ is called \emph{dominated} by $c(f) \in \R^{d}$ if $c(f) \leqslant c(f') .$ 
Accordingly, a feasible flow $f$ is called an \emph{efficient flow} if there exists no other feasible flow $f^{\prime}$ such that  $c ({f}^{\prime})  \leqslant c(f)$. 
The image $c(f)$ of an efficient flow $f$ is called a \emph{nondominated vector}. The set of efficient flows is denoted by $\mathcal{X}_E\subset \mathcal{X}$ and the set of nondominated vectors by $\mathcal{Y_N}\subset \mathcal{Y}$. Moreover, a feasible flow $f$ is called weakly efficient if there is no other flow $f^{\prime}$ such that $c({f}^{\prime}) < c(f)$\, where ``$<$'' denotes the strict componentwise order:
\[
c(f')  < c(f) \iff c^k(f') < c^k(f) \quad \forall k \in \{1,\dots , d\} 
\]

The polyhedron $\mathcal{P}=\conv(\mathcal{Y})+\R^d_{\geqq}$ is called the \emph{upper image} of \(\mathcal{Y}\).

Let $\|x\|_{1}$ denote the $1$-norm of $x \in \R^{d}$, i.\,e.,  $\|x\|_{1}\coloneqq \sum_{i=1}^{d}\left|x_{i}\right|$.
For the weighted sum method, we define the set of \emph{normalized weighting vectors} as the set $\Lambda_d=\{\lambda \in \R^{d}_{>} \colon \|\lambda\|_{1} = 1 \}$ or $\Lambda_d^0=\{\lambda \in \R^{d}_{\geq} \colon \|\lambda\|_{1} = 1 \}$ if weights equal to zero are included. Then, the \emph{weighted-sum linear program} w.\,r.\,t.\ $\lambda \in \Lambda_d$ or $\lambda \in \Lambda_d^0$ is defined as the parametric program 
\begin{equation}\tag{$P_{\lambda}$}\label{eq:P_l}
\begin{array}{rl}
\min &\lambda^\top  C\, f \\
\text{s.t.} & f \in \mathcal{X}
\end{array}
\end{equation} 

Note that the problem \eqref{eq:P_l} is as easy to solve as the associated single objective problem, as long as the encoding lengths of the components of $\lambda$ are not too large. (Consider that $\lambda$ is not part of the original input.) However, on a positive note, it is shown that these lengths can be bounded by $\mathcal{O}(\operatorname{poly}(n,m))$~\cite{Boekler2015}. If  $\lambda\in \Lambda_{d}^0$, every optimal solution to \eqref{eq:P_l} is a weakly efficient solution. Moreover, every optimal solution of \eqref{eq:P_l} is efficient, if \(\lambda\in\Lambda_{d}\) \cite{ehrgott2005multicriteria}. Considering a  MOMCF problem, the weighted sum scalarization achieves completeness when varying $\lambda$ within $\Lambda_p$~\cite{Iser74m} and any efficient solutions lie on the nondominated frontier, which is defined as the set $\{y \in \conv(\Y_N) \colon \operatorname{conv}(\Y_N) \cap (y -\R^p_{\geqq}) = \{y\}\}.$ 

Since the capacities on each arc are bounded, $\mathcal{X}$ is a polytope, and therefore, $\mathcal{Y}$ is also a polytope. We define the concepts of face, facet, maximal, and maximally nondominated faces similarly as in (\cite{EUSEBIO200968},~\cite{Steuer2008MultipleCO}). We call $F\subset \mathcal{Y}$ a \emph{face} of $\mathcal{Y}$ if there exists a supporting hyperplane $H$ such that $H \cap \mathcal{Y}=F.$ A face $F$ of $\mathcal{Y}$ is denoted as  \emph{$r$-face} if $F$ is of dimension $r$. Extreme points are 0-faces, and 1-faces are edges. A $r$-face $F\neq\mathcal{Y}$ is denoted as \emph{maximal face} or \emph{facet} of \(\mathcal{Y}\) if there does not exist another $s$-face $G$, such that $F \subset G$ and $r<s$. Note that the dimension of a facet $F$ 
 is $\dim(F)=\dim(\mathcal{Y})-1$. A face is called \emph{nondominated} if there exists a corresponding weight vector to the face with weights strictly greater than zero. A $r$-face $F$ is a \emph{maximally nondominated face} if there is no nondominated $s$-face $G$, such that $F \subset G$ and $r<s$. Let $F_\mathcal{X}$ be the set of the preimage of a maximally nondominated face $F_\mathcal{Y}$ of polytope $\mathcal{Y}$, i.\,e., all solutions whose image lies in $F_\mathcal{Y}$. We use the expression \emph{maximally efficient face} for $F_\mathcal{X}$. Even if we use the expression face for $F_\mathcal{X}$, it must not really be a face of the decision space since more than one feasible solution could have a vertex of $\mathcal{Y}$ as image, and they may lie in different faces of $\mathcal{X}$. The nondominated frontier equals the nondominated set of $\operatorname{conv} (\Y_N)$. It can be described as the set of the \emph{maximal nondominated faces} of the upper image $\mathcal{P}$~\cite{ehrgott2005multicriteria}.

However, in MOIMCF problems, nondominated points located in the interior of the upper image can occur, i.e., solutions whose image does not lie on the nondominated frontier and cannot be identified using a weighted-sum method. 
  These so-called \emph{unsupported points} typically outnumber these points that can be obtained with the weighted sum method~\cite{Visee1998}. As a result, determining the set of efficient solutions for MOIMCFP is significantly more challenging than for MOLP. Therefore, it is important to distinguish between different types of efficient solutions.

\begin{definition}\label{def:solution_types} The set of efficient solutions $\X_E$ can be distinguished in the following classes: 
\begin{enumerate}
\item \emph{Weakly supported efficient solutions} are efficient solutions that are optimal solutions of $P_{\lambda}$ for  $\lambda\in \Lambda_{p}^0$, i.e., an optimal solution to a weighted sum single-objective problem with weights strictly or equal to zero. Their images in the objective space are weakly \emph{supported nondominated points}; we use the notations $\X_{wS}$ and $\Y_{wS}$. All weakly supported nondominated vectors are located on the boundary of the \emph{upper image} $\mathcal{P}\coloneqq  \operatorname{conv}(\Y)+ \R_\geqq^p.$

\item An efficient solution is called \emph{supported efficient solution} if it is an optimal solution to the weighted sum scalarization $P_{\lambda}$ for  $\lambda\in \Lambda_{d}$, i.e., an optimal solution to a single objective weighted-sum problem where the weights are strictly positive. Its image is called \emph{supported nondominated vector}; we use the notations $\X_{S}$ and $\Y_{S}$ and denote the cardinality of $\mathcal{X}_S$ by $|\mathcal{X}_S| = \mathcal{S}$. A supported nondominated vector is located on the \emph{nondominated frontier}, i.e., located on the union of the \emph{maximal nondominated faces}. 

\item \emph{Extreme supported solutions}  are those solutions whose image lies on the vertex set of the upper image. Their image is called an \emph{extreme supported nondominated vector}.

\item \emph{Unsupported efficient solutions} are efficient solutions that are not optimal solutions of  $P_{\lambda}$ for any  $\lambda\in \Lambda_{p}^0$. \emph{unsupported nondominated vectors} lie in the interior of the upper image. 
\end{enumerate}
\end{definition}

When solving multi-objective problems, they can be solved concerning different subsets of interest, and it is important to specify what class of set of efficient solutions is searched and be aware of the (partial or complete) enumeration of equivalent solutions. For example, the complete set of all efficient solutions $\X_E$ can be searched. However, if we speak about the task to determine all nondominated vectors $\mathcal{Y_N}$,  we want to determine a minimal complete set of $\X'\subseteq \X_E$.  If a subset $\X'\subseteq \X$ fulfills $f(\X')=\mathcal{Y_N}$, i.e., for each $y\in \mathcal{Y_N}$ there exists at least one $x'\in \X'$ such that $f(x)=y$, $\X'$ is called a \emph{complete set} of efficient solutions. A complete set is called \emph{minimal complete set} if it contains no  equivalent solutions. The set of all efficient solutions $\X_E$ is then the \emph{maximal complete set} of efficient solutions.

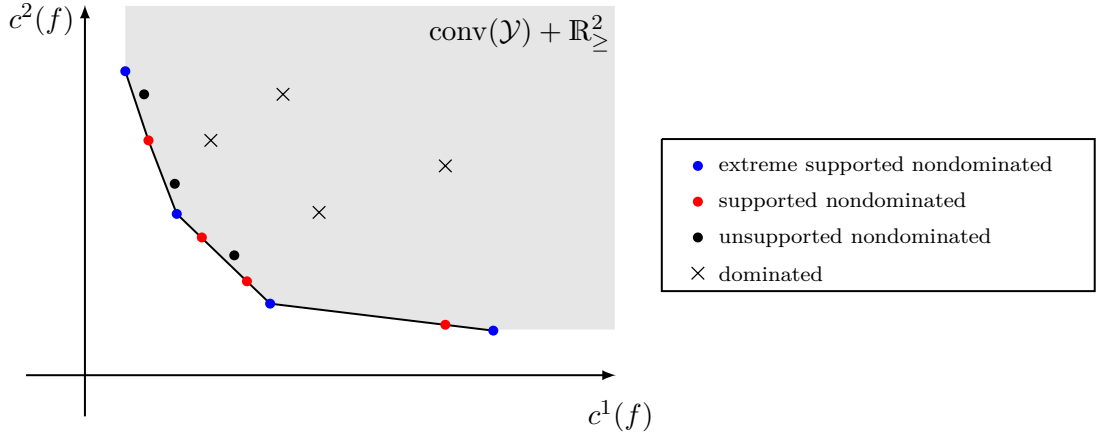
\begin{figure}[t]
\centering
\tikzset{every picture/.style={line width=0.75pt}} %set default line width to 0.75pt        
\tikzstyle{vertex}=[circle,fill=blue,draw=blue,minimum size=3pt,inner sep=0]
\tikzstyle{vertex2}=[circle,fill=red,draw=red,minimum size=3pt,inner sep = 0]
\tikzstyle{vertex3}=[circle,fill=black,draw=black,minimum size=3pt,inner sep = 0]

\newcommand{\Cross}{$\mathbin{\tikz [x=1ex,y=1ex,line width=.1ex, black] \draw (0,0) -- (1,1) (0,1) -- (1,0);}$}%
\newcommand{\Crossa}{$\mathbin{\tikz [x=1ex,y=1ex,line width=.1ex, black] \draw (0,0) -- (1,1) (0,1) -- (1,0);}$}%
\newcommand{\Crossplus}{$\mathbin{\tikz [x=1ex,y=1ex,line width=.1ex, black] \draw (0,0.5) -- (1,0.5) (0.5,0) -- (0.5,1);}$}%
\newcommand{\Triangle}{$\mathbin{\tikz [x=1ex,y=1ex,line width=.1ex, blue] \draw (0,0) -- (0.5,1) -- (1,1) -- (0,0);}$}%
\begin{tikzpicture}[x=0.75pt,y=0.75pt,yscale=-1,xscale=1,scale=0.9]
%uncomment if require: \path (0,353); %set diagram left start at 0, and has height of 353

%Shape: Axis 2D 
\draw[-latex] (127.5,257) -- (454.5,257);
\draw[latex-] (160.2,50) -- (160.2,280)  ;

%Shaded Area:
\draw [draw= white, fill= black!10  ,fill opacity=1 ]   (182.58,87.08) -- (195.42,125.75) -- (211.11,166.83) -- (262.87,217) -- (386.58,232.08) -- (454.5,232.08) -- (454.5,50) -- (182.58,50);

%Straight Lines [id:da8047586419045308] 
\draw [fill={rgb, 255:red, 255; green, 255; blue, 255 }  ,fill opacity=0 ]   (182.58,87.08) -- (195.42,125.75) -- (211.11,166.83) -- (262.87,217) -- (386.58,232.08) ;
%Shape: Nodes
\node[vertex] at (182.58,87.08) {};
\node[vertex] at (211.11,166.83) {};
\node[vertex] at (262.87,217) {};
\node[vertex] at (386.58,232.08) {};

%\node at (386.58,232.08) {\Triangle};%

 %\node at (210,150) {\Crossplus};
  %\node at (193,100) {\Crossplus};
   %\node at (243,190) {\Crossplus};

\node[vertex2] at (195.42,125.75) {};
\node[vertex2] at (250,204.5) {};
\node[vertex2] at (225,180) {};
\node[vertex2] at (360,228.8) {};
\node[vertex3] at (210,150) {};
\node[vertex3] at (243,190) {};
\node[vertex3] at (193,100) {};

%Cross Nodes:
%\node at (210,150) {\Cross};
%\node at (243,190) {\Cross};
%\node at (193,100) {\Cross};
%\node at (300,210) {\Cross};

\node at (290,166) {\Crossa};
\node at (230,125.75) {\Crossa};
\node at (270,100) {\Crossa};
\node at (360,140) {\Crossa};

% Text Node
\draw (112,43.4) node [anchor=north west]    {$c^{2}(f)$};
% Text Node
\draw (433,265) node [anchor=north west]    {$c^{1}(f)$};

\draw (345,50) node [anchor=north west] {$\conv(\mathcal{Y})+\R^2_{\geq}$};

%Draw legend
\begin{scope}[shift={(320,-180)}]
\draw (160,305) -- (400,305) -- (400,390) -- (160,390) -- (160,305);
\node[vertex] at (180,320) {};
\draw (190,320) node [anchor= west][inner sep=0.75pt]   [align=left] {\scriptsize extreme supported nondominated};
\node[vertex2] at (180,340) {};
\draw (190,340) node [anchor= west][inner sep=0.75pt]   [align=left] {\scriptsize supported nondominated};
\node[vertex3] at (180,360) {};
\draw (190,360) node [anchor= west][inner sep=0.75pt]   [align=left] {\scriptsize unsupported nondominated};
\node at (180,380) {\Crossa};
\draw (190,380) node [anchor= west][inner sep=0.75pt]   [align=left] {\scriptsize dominated};
\end{scope}

\end{tikzpicture}
\caption{ Illustration of the upper image $\mathcal{P}=\conv(\mathcal{Y})+ \R^2_{\geqq}$ and the different solution types.}\label{Image:DiffSolutions}\end{figure}

\Cref{Image:DiffSolutions} illustrates supported extreme, supported, and unsupported nondominated points, as well as the upper image in the bi-objective case.

Due to the total unimodularity of MOMCF, there exists for each extreme supported nondominated point of MOMCF an integer flow mapping to it, if $u$ and $b$ are integral (see~\cite{Ahuja1993}). In other words, the sets of extreme supported nondominated points of MOMCF and MOIMCF, and thus the respective upper images, coincide.
In the remainder of this paper, only integer flows are considered, and from now on, flow always refers to an integer flow.  

The unsupported nondominated points typically outnumber the supported ones in MOCO problems and are more difficult to obtain. However, as shown in~\cite{Serpil2024}, the supported or extreme supported point set yields good representations for MOCO problems. The distinction between supported and weakly supported solutions is relatively new, first presented in~\cite{koenen2025supportednessmultiobjectivecombinatorialoptimization}, to provide a consistent definition of supportedness in MOCO problems. \Cref{Ex:weakly} illustrates an instance of the MOIMCF problem where weakly supported nondominated points exist, and the supported nondominated vectors form a proper subset of the weakly supported solutions.

\begin{example}\label{Ex:weakly}
    The MOIMCF problem shown in~\Cref{example1} has three extreme nondominated vectors $y^1=(8,16,6)^\top , y^2=(12,12,6)^\top ,y^3=(16,8,10)^\top$ and in addition the following nondominated vectors $s^1=(9,15,7)^\top , s^2=(10,14,8)^\top ,s^3=(11,13,9)^\top  ,s^4=(13,11,7)^\top , s^5=(14,10,8)^\top ,s^6=(15,9,9)^\top $ and as well one dominated vector $d^1 = (12,12,10)^\top $. The different vectors in objective space and their convex hull are illustrated in~\Cref{obj-ex1}. All vectors $s^i$ with $i\in \{1,\ldots, 6\}$  are weakly supported nondominated vectors since their preimages are optimal solutions of the corresponding weighted sum problem $P_\lambda$ with $\lambda= (0.5 , 0.5, 0)^\top $. The vectors $s^4,s^5,s^6$ are also supported nondominated vectors since they are optimal solutions of $(P_{\lambda^2})$ with $\lambda^2=(0.25,0.5,0.25)^\top $. However, $s^1,s^2,s^3$ are not optimal for any weighted sum problem $P_\lambda$ with $\lambda\in \Lambda_d$.  Hence, the set of weakly supported nondominated vectors is $\mathcal{Y}_wS=\{y^1,y^2,y^3,s^1,s^2,s^3,s^4,s^5,s^6 \}$. While the set of supported nondominated vectors is $\mathcal{Y}_{S}=\{y^1,y^2,y^3,s^4,s^5,s^6 \}$. Thus, in this example, the set of supported nondominated vectors is a proper subset of the set of weakly supported nondominated vectors $\mathcal{Y}_{S} \subset \mathcal{Y}_{wS}$. It is easy to see that there do not exist weights $\lambda\in \Lambda_d$ such that $s_1$, $s_2$, or $s_3$ are optimal solutions for $P_\lambda$ using the \emph{weight space decomposition} which will be formally introduced in~\Cref{sec: multi}. \Cref{fig:wsd} shows the weight space decomposition of the \emph{upper image} of the MOIMCF instance given in \Cref{example1}.
\end{example}

\tikzstyle{vertex}=[circle,fill=black,draw=black,minimum size=8pt,inner sep=0]%
\tikzstyle{edge} = [draw,thick,->]%
\tikzstyle{selected vertex} = [vertex, fill=red!24]%
\tikzstyle{selected edge} = [draw,line width=5pt,-,yellow!50]%
\usetikzlibrary{arrows,automata}%
\pgfdeclarelayer{background}%
\pgfsetlayers{background,main}%

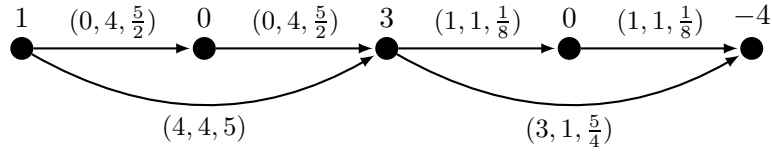
\begin{figure} 
\centering
	\centering
		\begin{tikzpicture}[scale=1.2,-latex,shorten >=1pt,auto, semithick]%
					\foreach \pos/\name/\label in {{(-2,0)/a/1}, {(-0,0)/b/0}, {(2,0)/c/3},%
						{(4,0)/d/0}, {(6,0)/e/-4}}%
					\node[vertex] (\name) [label=$\label$] at \pos {};%
					\path (a) edge[-latex,thick]   	node[] {\small$(0,4,\frac{5}{2})$} (b)%
					(b) edge[-latex,thick]   	node[] {\small$(0,4,\frac{5}{2})$} (c)%
					(c) edge[-latex,thick]   	node[] {\small$(1,1,\frac{1}{8})$} (d)%
					(d) edge[-latex,thick]   	node[] {\small$(1,1,\frac{1}{8})$} (e)%
					(a) edge[-latex,thick,bend right]   	node[sloped,below] {\small$(4,4,5)$} (c)%
					(c) edge[-latex,thick,bend right]   	node[sloped,below] {\small$(3,1,\frac{5}{4})$} (e);
					\end{tikzpicture}%
		\caption{The Graph corresponds to a tri-objective MOIMCF problem. Thereby,  $l_{ij}=0$ and  $u_{ij}=4$ for all arcs $(i,j)\in A$. The arcs are labeled with their cost coefficients $(c^1_{ij},c^2_{ij},c^3_{ij})$ and the nodes are labeled with their supply/demand $b_i$.}\label{example1}
\end{figure}

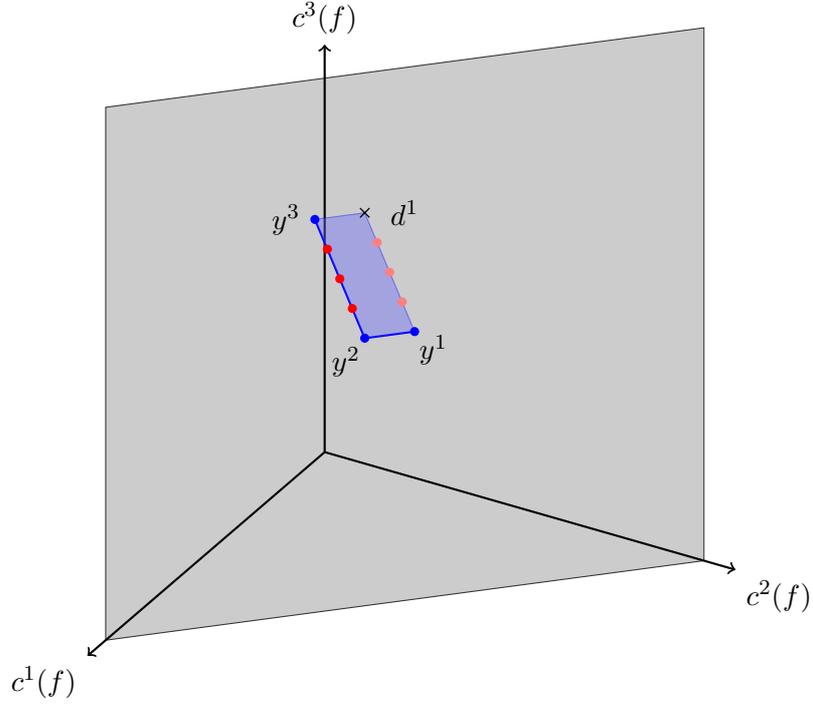
\begin{figure}
\centering
%\tdplotsetmaincoords{70}{110}
\tikzstyle{vertex}=[circle,fill=SteelBlue2,draw=SteelBlue2,minimum size=3pt,inner sep=0]
\tikzstyle{vertex2}=[circle,fill=DeepPink2,draw=DeepPink2,minimum size=3pt,inner sep=0]
\tikzstyle{vertex3}=[circle,fill=DeepPink2!50,draw=DeepPink2!50,minimum size=3pt,inner sep = 0]
\tikzstyle{vertex4}=[circle,fill=black,draw=black,minimum size=3pt,inner sep=0]
\newcommand{\Crossa}{$\mathbin{\tikz [x=1ex,y=1ex,line width=.1ex, black] \draw (0,0) -- (0.75,0.75) (0,0.75) -- (0.75,0);}$}%
\begin{tikzpicture}[scale=0.14,tdplot_main_coords, bluenode/.style={color=SteelBlue2,circle,fill,minimum size=2.5pt,inner sep=0pt}, blacknode/.style={circle,fill,minimum size=2.5pt,inner sep=0pt}
	every label/.style={font=\footnotesize}]%scale=3,tdplot_main_coords]
   
    \draw[thick,->] (0,0,0) -- (52,0,0) node[anchor=north east]{$c^1(f)$};
    \draw[thick,->] (0,0,0) -- (0,52,0) node[anchor=north west]{$c^2(f)$};
    \draw[thick,->] (0,0,0) -- (0,0,52) node[anchor=south]{$c^3(f)$};

     \filldraw[
        draw=black!15,%
        fill=black!15,%
        opacity= 0.8,
    ]          (24,0,0)
            -- (0,24,0)
            -- (0,24,55)
            -- (24,0,55)
            -- cycle;

	\draw[thick,SteelBlue4] (16,8,40) -- (12,12,24);
	\draw[thick,SteelBlue4] (8,16,24) -- (12,12,24);
	\filldraw[
        draw=SteelBlue4,%
        fill=SteelBlue4,%
        opacity=0.4,
    ]          (12,12,40)
            -- (16,8,40)
            -- (12,12,24)
            -- (8,16,24)
            -- cycle;
	\node at (12,12,40) [label=right:{$d^1$}]  (d1) {\Crossa};
	\node[vertex] at (8,16,24) [label={[xshift=7pt, yshift=-18pt]$y^1$}]  (y1) {};
	\node[vertex] at (12,12,24) [label={[xshift=-7pt, yshift=-20pt]$y^2$}]  (y2) {};
	\node[vertex] at (16,8,40) [label=left:{$y^3$}]  (y3) {};
        \node[vertex3] at (9 ,15 ,28) [label=right:{$s^1$}] {};
        \node[vertex3] at (10 ,14 ,32)  [label=right:{$s^2$}]{};
        \node[vertex3] at (11 ,13 ,36)  [label=right:{$s^3$}]{};
        \node[vertex2] at (15 ,9 ,36)  [label=left:{$s^6$}]{};
        \node[vertex2] at (14 ,10 ,32)  [label=left:{$s^5$}]{};
        \node[vertex2] at (13 ,11 ,28)  [label=left:{$s^4$}]{};

        %\node[vertex4] at (24 ,0 ,0)  [label=left:{}]{};
        %\node[vertex4] at (24 ,0 ,55)  [label=left:{}]{};
        
        %\node[vertex4] at (0 ,24 ,0)  [label=right:{}]{};
        %\node[vertex4] at (0 ,24 ,55)  [label=right:{}]{};

        \draw (24,0,0) -- (24,0, 55); 
        \draw (0,24,0) -- (0,24, 55); 
        \draw (24,0,0) -- (0,24,0);

\end{tikzpicture}
\caption{The figure shows the convex hull $\conv(\mathcal{Y})$ of the given problem in~\Cref{example1} in blue, all its integer vectors and the hyperplane $h=\{y\in \mathcal{Y} \colon 0.5\,y_1+0.5\,y_2=12\}$ in gray. The red dots on the edge between $(y^2,y^3)$ are the vectors $s^4,s^5$, and $s^6$. The light red dots on the edge between the edge $(y^1,d^1)$ are the vectors $s^1,s^2$, and $s^3$, which would be dominated in the non-integer case by the points on the edge $(y^1,y^2)$. The maximally nondominated faces are the edges $(y^1,y^2)$ and $(y^2,y^3)$. Note that $s^1$,$s^2$, and $s^3$ would lie on $\conv(\mathcal{Y})$  but not in any of the maximally nondominated faces. The vectors $s^1,s^2$, and $s^3$ would also lie on the boundary of the upper image. }
\label{obj-ex1}
\end{figure}

The clear distinction between the sets of supported nondominated and weakly supported nondominated solutions is also necessary as the corresponding problems differ in their output time complexity. 
 Eus\'{e}bio and Figueria~\cite{EUSEBIO200968} rely on a definition, which is equivalent to what we define as weakly supported nondominated vectors.  However, their proposed algorithm computes only supported efficient solutions whose images lie in the maximally nondominated faces.
In this paper, an output-polynomial algorithm is developed to determine all supported efficient solutions to the MOIMCF problem  with a fixed number of objectives. However, unless $\mathbf{P} = \mathbf{N P}$, the next chapter proves that an output-polynomial algorithm to determine all weakly supported nondominated vectors (or all weakly supported efficient solutions) for a MOIMCF problem with a fixed number of more than two objectives can be excluded. This is proven by showing that an output-polynomial algorithm cannot exist for the multi-objective $s$-$t$-path problem (MOSP) with more than two objectives.

\section{Output-sensitivity Analysis for the Weakly Supported Efficient Solutions}

This chapter starts by formally introducing the theory of \emph{output-sensitive} complexity of enumeration problems. 
For a comprehensive introduction, see~\cite{JOHNSON1988119}. The concepts discussed here are based on~\cite{boeklerthesis}. 

\subsection{An Introduction to Output-Sensitive Complexity}

\begin{definition}[\cite{Boekler2017}]
     An enumeration problem is a pair $(I, C)$ such that 
\begin{itemize}
    \item[1.] $I \subseteq \Sigma^*$ is the set of \emph{instances} for some fixed alphabet $\Sigma$,
    \item[2.] $C\colon I \rightarrow 2^{\Sigma^*}$ maps each instance $x \in I$ to its \emph{configurations} $C(x)$, and
    \item[3.] the encoding length $|s|$ for $s \in C(x)$ for $x\in I$ is in $\operatorname{poly}(|x|)$,
    \end{itemize}
    where $\Sigma^*$ can be interpreted as the set of all finite strings over $\{0,1\}$. 
\end{definition}
We assume that $I$ is decidable in polynomial time and that $C$ is computable.
Then an \emph{enumeration algorithm} for an enumeration problem $E=(I, C)$ is a \emph{random access machine} that on input $x \in I$ outputs $c\in C(x)$ exactly once, and
 on every input terminates after a finite number of steps.
An enumeration problem is called \emph{intractable} if the cardinality of the configuration set is exponential in the size of an input instance.

MOIMCF problems, like any MOCO problem, can be considered as enumeration problems. This involves enumerating all nondominated points or subsets of the nondominated point set with respect to other solution concepts.

Optimization problems are often classified based on whether they can be solved in polynomial or exponential time. However, many MOCO problems, as the MOIMCF problem, are known to be intractable, i.e., the cardinality of the set of efficient solutions or the set of nondominated points grows exponentially with the input size. Even when finding a single solution may be easy, or the corresponding canonical decision problem may be solvable in polynomial time, enumerating all desired solutions of such an intractable problem in polynomial time is impossible.

Thus, for such problems, one is interested in whether the enumeration problem can be solved in output-polynomial time, which can be formally defined by:

\begin{definition}[\cite{Boekler2017}]
    An enumeration algorithm for an enumeration problem $ E=(I, C)$ is said to run in \emph{output-polynomial time} \emph{(is output-sensitive)} if its running time is in $\operatorname{poly}(|x|,|C(x)|)$ for $x\in I$.
\end{definition}

Let $C$ be the configuration set of all supported efficient solutions and $C^*$ the configuration set of all nondominated vectors. An output-polynomial time algorithm $E=(I,C)$, which determines all supported efficient solutions exactly once, would also yield all nondominated vectors. However, it might not be output-polynomial w.r.t.\ this task, as there might be an exponential number of efficient solutions mapping to a small number of nondominated points. Thus, the algorithm for $E^*=(I,C^*)$  may output the elements (i.e., the nondominated vectors) more than once or even exponentially many times. Consider, for example, a directed graph with $\{1,\dots,n\}$ transshipment nodes, a node $s$ and $t$ with flow balance $n$ and $-n$, respectively. The graph contains the arcs $(s,i)$ and $(i,t)$ for all  $i \in\{1,\dots,n\}$ with upper capacity $n$. The cost of all arcs is equal to \((1,1)^\top\in\R^2\). Then we have $\binom{2n-1}{n}$ supported efficient solutions, but all map to the same extreme nondominated point. 

A finished decision problem $E^{\mathrm{Fin}}$ for an enumeration problem $E=(I, C)$ is defined as follows: Given an instance $x \in I$ of the enumeration problem and a subset $M \subseteq C(x)$ of the configuration set, the goal is to decide if $M=C(x)$, i.e.,  to determine if all configurations have been found.  
 If the enumeration problem $E$ can be solved in output-polynomial time, then $E^{\mathrm{Fin}}\in \mathbf{P}$~\cite{Lawer1972}.

\subsection{Weakly Supported Efficient Solutions}

\begin{thm}[\cite{koenen2025supportednessmultiobjectivecombinatorialoptimization}]\label{thm:d+1}
    The determination of all weakly-supported solutions for a MOCO problem with $d+1$ objectives is as hard as the determination of all nondominated points for a MOCO problem with $d$ objectives. 
\end{thm}

It is well known that the multi-objective shortest path problem, a special case of the minimum cost flow problem, is intractable and there could not be an output-polynomial algorithm to determine all nondominated points, even in the bi-objective case. 

 Given a directed graph $D = (V, A)$, a cost function $c\colon A \rightarrow \mathbb{R}^p$ that assigns a cost vector to each arc, a source node \(s \in V\), and a target node \(t \in V\), the  (multi-objective) \emph{shortest path} problem is defined as:

\begin{equation}\label{eq:MOSP}
        \tag{MOSP}
     \min_{P\in\mathcal{P}_{s,t}}\left\{ c(P)=\left(c_1(P),\ldots,c_p(P)\right)^{\top} = Cx\right\},     
\end{equation}
    where \(\mathcal{P}_{s,t}\) is the set of all paths from \(s\) to \(t\) in \(G\),  $C \in  \R^{p\times n}$ is the cost matrix containing the rows $c^k$ of coefficients of the $p$ linear objective functions and $c_i(P) = \sum_{e\in P} c^i_e$ for $i\in \{1,\ldots,p\}$.

\begin{lem}[\cite{Boekler2017}]\label{lem:mosp}
 There is no output-polynomial algorithm for determining all nondominated vectors for the~\ref{eq:MOSP} unless $\mathbf{P} = \mathbf{N P}$.
\end{lem}

The determination of all nondominated vectors of a \ref{eq:MOSP}, as well as the determination of all efficient solutions of a \ref{eq:MOSP} can be formulated as an enumeration problem~\cite{Boekler2017}. We denote the finished decision problem for the determination of all nondominated vectors of a \ref{eq:MOSP} as \ref{eq:MOSP}$^{\mathrm{Fin}}_{\mathcal{Y}}$ and the determination of all efficient solutions of a \ref{eq:MOSP} as  \ref{eq:MOSP}$^{\mathrm{Fin}}_{\mathcal{X}}$, respectively. 
The proof of~\Cref{lem:mosp} is given by showing that the finished decision variant \ref{eq:MOSP}$^{\mathrm{Fin}}_{\mathcal{Y}}$ of the \ref{eq:MOSP} is $\mathbf{co}$-$\mathbf{NP}$-hard, by a reduction  of the complement of the Knapsack problem:
\begin{equation}\label{eq:KP}\tag{KP}
\Bigl\{\bigl(c^1, c^2, k_1, k_2\bigr) \colon c^{1^\top } x \leq k_1, c^{2^\top } x \geq k_2, x \in\{0,1\}^n\Bigr\}.
\end{equation}
This problem is  NP-complete~\cite{KelPfePis04}.

The same reduction that is used in the proof of~\Cref{lem:mosp} given in~\cite{Boekler2017} can be extended to show that the \ref{eq:MOSP}$^{\mathrm{Fin}}_{\mathcal{X}}$ is also $\mathbf{co}$-$\mathbf{NP}$-hard and thus there cannot exist an output-polynomial time algorithm to determine all efficient solutions for the \ref{eq:MOSP}. However, we have to adjust the costs of some weights. In~\Cref{fig:proof-fritz}, an example of the reduction is given,  similar to the one used in the proof in~\cite{Boekler2017} to show the following Lemma. For the sake of simplicity, we do not give a formal proof here. For more details, we refer to~\cite{Boekler2017}.

 \tikzstyle{svertex}=[circle,fill=black,draw=black,minimum size=5pt,inner sep=0]%
\tikzstyle{wvertex}=[circle,fill=white,draw=white,minimum size=8pt,inner sep=0]
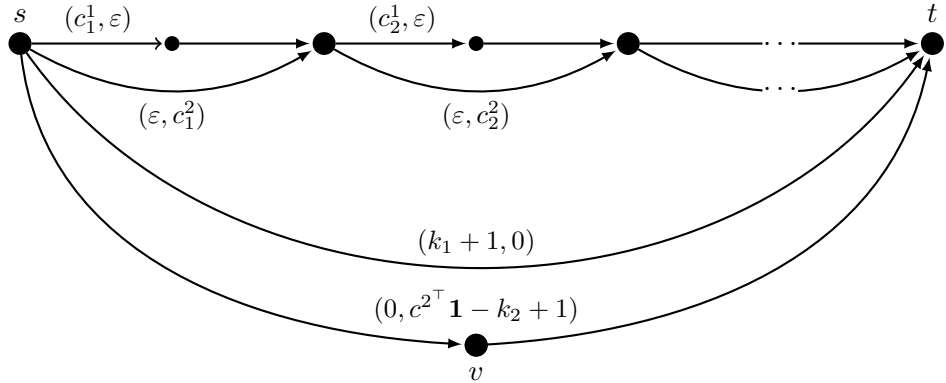
\begin{figure} 
	\centering
	\centering
	\begin{tikzpicture}[ ->,shorten >=1pt,auto,node distance=2.8cm,%
		semithick]%
		\foreach \stil/\pos/\name/\label in {{vertex/(-2,0)/a/s}, {svertex/(-0,0)/b/}, {vertex/(2,0)/c/},%
			{svertex/(4,0)/d/}, {vertex/(6,0)/e/}, {wvertex/(8,0)/f/}, {vertex/(10,0)/g/t},{vertex/(4,-4)/v/$\small$(0, c^{2^\top} \mathbf{1}-k_2+1)}}%
		\node[\stil] (\name) [label=$\label$] at \pos {};%

		\path (a) edge[->,thick]  node[sloped] {\small$(c^1_1,\varepsilon)$}  	 (b)%
		(b) edge[-latex,thick]   	 (c)%
		(c) edge[-latex,thick]   node[sloped] {\small$(c^1_2,\varepsilon)$}  (d)%
		(d) edge[-latex,thick]   	 (e)%
            (e) edge[-,thick]   	 (f)%
            (f) edge[-latex,thick]   	 (g)%
		(a) edge[-latex,thick,bend right]  node[sloped,below] {\small$(\varepsilon,c^2_1)$} (c)
		(c) edge[-latex,thick,bend right]  node[sloped,below] {\small$(\varepsilon,c^2_2)$} (e)
            (a) edge[-latex,thick,bend right,out=310]   (v)
            (e) edge[-latex,thick,bend right]   (g)
            (v) edge[-latex,thick,bend right,in=225]   (g)
            (a) edge[-latex,thick,bend right, in=235,out=305]  node[sloped,above] {\small$(k_1+1,0)$} (g)
		%(e) edge[->,thick,bend left]  node[sloped,below] {$(\ell,0,1,-1)$} (c);
            ;

            \node[wvertex] (h) at (8,-0.6) {$\ldots$}; 
            \node[wvertex] (h1) at (8,0) {$\ldots$};  
            \node[vertex,label=below:$v$]  (v1) at (4,-4) {};  
	\end{tikzpicture}%
	\caption{Showing the reduction in the proof of ~\Cref{lem:mosp} prestented by Bökler et al.~\cite{Boekler2017} with modified costs in order to proof~\Cref{lem:mosp_x}. Here $\varepsilon= 1/(n+1)$ and arcs with no label have cost $\mathbf{0}$. Determining another efficient $s$-$t$ path to the two with cost equal to $ (k_1+1,0) $ and $ (0, c^{2^\top} \mathbf{1}-k_2+1) $ would be an instance of the \eqref{eq:KP} problem. }\label{fig:proof-fritz}
\end{figure}

\begin{lem}[]\label{lem:mosp_x}
 There is no output-polynomial time algorithm to determine all efficient solutions for the~\ref{eq:MOSP} problem unless  $\mathbf{P} = \mathbf{N P}$.
\end{lem}

The multi-objective shortest path problem~(\ref{eq:MOSP}) is a special case of the multi-objective minimum cost flow problem~\cite{Ahuja1993}.  The $s$-$t$ path problem can be transformed into a minimum cost flow problem by setting all arc capacities to one and sending one unit of flow from $s$ to $t$, with $b_s=1$, $b_t=-1$, and $b_v=0$ for all $v \in V\backslash\{s,t\}$. In this modified network, an optimal flow corresponds to the shortest path from $s$ to $t$. Consequently, the finished or canonical decision problem of the multi-objective minimum cost flow problem is $\mathbf{NP}$-hard, and the enumeration problem is intractable even in the case of two objectives. 
Showing that another nondominated vector or efficient flow exists would solve the complement of the knapsack problem.

Thus, due to~\Cref{lem:mosp},~\Cref{lem:mosp_x} and~\Cref{thm:d+1}, it holds that there is no output-polynomial algorithm to determine all \weakly supported nondominated vectors or weakly supported efficient solutions for the MOIMCF problem with a fixed number of $d \geq 3$ objectives.

\begin{thm}\label{thm:moimcf}
Unless $\mathbf{P} = \mathbf{N P}$, there is no output-polynomial algorithm  to determine all \emph{weakly} supported nondominated vectors (or weakly supported efficient solutions) for the MOIMCF problem with a fixed number of $d \geq 3$ objectives.
\end{thm}

\section{Finding all \strictly Supported 
Efficient Flows} 
\label{3}
In this section, output-polynomial time algorithms are derived
that determine all supported efficient flows for MOIMCF. It is based on the determination of all optimal flows for a sequence of single-objective parametric network flow problems, each corresponding to a maximally nondominated face. The approach consists of two phases and relies on the following widely known fact for any integer supported flow, see, e.g., \cite{gal77}.

\begin{thm}\label{thm:supp_face}
A flow $f$ is contained in $F_\mathcal{X}$, i.\,e., its image of $c(f)$ lies on a maximally nondominated face $F_\mathcal{Y}$ w.\,r.\,t.\ an associated weight vector $\lambda \in \Lambda_d$, if $f$ is an optimal solution to the parametric network flow program  \eqref{eq:P_l}.
\end{thm}

Any image of a \strictly supported flow must lie in at least one maximally nondominated face, and any integer point in a maximally nondominated face corresponds to a \strictly supported integer flow. 
Assuming that for each maximally nondominated face $F_i\in\{F_1,\dots, F_t\}$ one optimal solution $f^i$ and the corresponding weighting vectors $\lambda^i$ are given the problem of determining all \strictly supported flows reduces to determining all optimal flows for each parametric single-objective problem  $(P_{\lambda^i})$. These optimal solutions can be determined by using the algorithm for determining all optimum flows for single-objective minimum cost flow problems presented in~\cite{KONEN2022333}, which we refer to as the all optimum flow (AOF) algorithm. 
 The AOF algorithm successively searches for so-called \emph{proper zero-cost} cycles efficiently by using a modified depth-first search technique. 

\begin{thm}[\cite{KONEN2022333}]
Given an initial optimal integer flow $f$,  we can determine all optimal integer flows in $\mathcal{O}(\mathcal{F}(m+n)+m\, n)$ time for a single-objective minimum cost flow problem, where $\mathcal{F}$ is the number of all optimal integer flows. 
\end{thm}

Therefore, we divide the approach into two phases: In phase one, we determine all extreme nondominated vectors and one weighting vector for each maximally nondominated face. In phase two, we  apply the AOF algorithm to the corresponding weighted-sum program for each maximally nondominated face. 

The extreme nondominated vectors and the weighting vectors for each maximally nondominated face can be determined much easier in the case of two objectives. Thus, we start by deriving an algorithm for BOIMCF problems and consider the general case of MOIMCF problems afterward.

\subsection{Bi-Objective Minimum Cost Flow Problem}

In the bi-objective case, the set of supported flows is equal to the set of weakly supported flows since every weakly nondominated face contains exactly one nondominated vector (namely an extreme point), which dominates the complete face.

First, we determine all $N$ extreme points and precisely one corresponding extreme flow by using the enhanced parametric programming approach in  $\mathcal{O}(M+ Nn(m+n\log n))$ time~\cite{raith17}, where $M$ denotes the time required to solve a given single-objective minimum cost flow problem. Also, the algorithm stores one extreme flow for each extreme nondominated point. 

In the bi-objective case, every maximally nondominated face $F_\mathcal{Y}$ of $\conv(\mathcal{Y})$ is a line segment connecting two adjacent extreme supported points if there is more than one nondominated point (\(|\mathcal{Y}_N|>1\)). A maximally nondominated face can only have dimension zero if there is only one extreme nondominated point, which implies that there is only one nondominated point (or, in other words, the ideal point is feasible).

In the following, we will derive a procedure to determine the complete set of all supported efficient flows. For that, we will determine all supported flows whose images lie on the maximally nondominated edges. 
Let $y^1,\dots,y^N$ be the extreme supported points obtained by the enhanced parametric programming approach \cite{raith17} and let $f^1,\ldots, f^N$ be a set of corresponding extreme supported flows each mapping to one extreme supported point. Moreover, we sort the set of extreme supported points and flows $\{y^i=(c^1(f^i),c^2(f^i)),f^i\colon i\in\{1,\ldots,N\}\}$ by non-decreasing values of $c^1$.
For each pair of consecutive extreme points $y^i$ and $y^{i+1}$, we determine the weighting vector \(\lambda^i\in\Lambda\) that corresponds to the normal of the maximally nondominated facet \(F_i\) connecting the extreme points $y^i$ and $y^{i+1}$:

\begin{align*}
	\lambda^i \coloneqq 
	\begin{pmatrix}
		c^2(f^i)-c^2(f^{i+1})\\
		c^1(f^{i+1})-c^1(f^i)
	\end{pmatrix}
	% \frac{1}{\sqrt{\lambda_1^2+\lambda_2^2}} \,\bigl( \lambda_1, \lambda_2\bigr)^\top
	% \qquad
	% \text{where}\quad & 
	% \lambda_1^i \coloneqq c^2(f^i)-c^2(f^{i+1}) \\[-1ex]
	% & \lambda_2^i \coloneqq c^1(f^{i+1})-c^1(f^i). 
\end{align*}

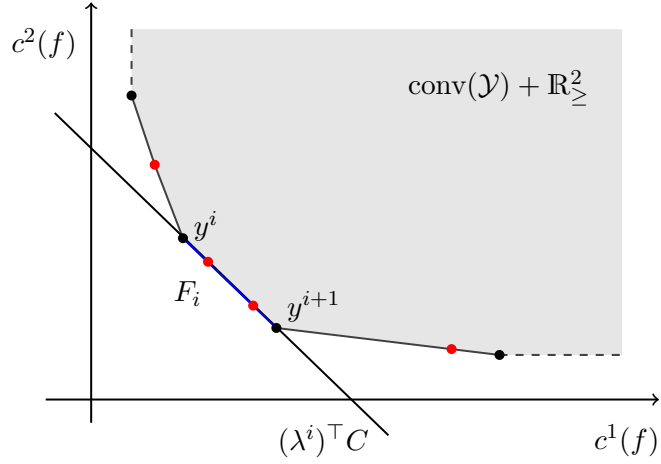
\begin{figure}[htb]
\centering
\tikzset{every picture/.style={line width=0.75pt}} %set default line width to 0.75pt        
\tikzstyle{vertex}=[circle,fill=black,draw=black,minimum size=3pt,inner sep=0]
\tikzstyle{vertex2}=[circle,fill=red,draw=red,minimum size=3pt,inner sep = 0]
\newcommand{\Cross}{$\mathbin{\tikz [x=1ex,y=1ex,line width=.1ex, black] \draw (0,0) -- (1,1) (0,1) -- (1,0);}$}%
\newcommand{\Crossa}{$\mathbin{\tikz [x=1ex,y=1ex,line width=.1ex, blue] \draw (0,0) -- (1,1) (0,1) -- (1,0);}$}%
\begin{tikzpicture}[x=0.75pt,y=0.75pt,yscale=-1,xscale=1,scale =0.9]
%uncomment if require: \path (0,353); %set diagram left start at 0, and has height of 353

%Shape: Axis 2D 
\draw [->] (135,257) -- (475,257);
\draw [->] (160.2,270) -- (160.2,35);
%(160.2,50) -- (160.2,280) (447.5,252) -- (454.5,257) -- (447.5,262) (155.2,57) -- (160.2,50) -- (165.2,57)  ;

%Shaded Area:
\fill [ fill= black!10  ,fill opacity=1 ]   (182.58,87.08) -- (195.42,125.75) -- (211.11,166.83) -- (262.87,217) -- (386.58,232.08) -- (454.5,232.08) -- (454.5,50) -- (182.58,50);

%Straight Lines [id:da8047586419045308] 
\draw [color = darkgray, fill=gray  ,fill opacity=0 ]   (182.58,87.08) -- (195.42,125.75) -- (211.11,166.83) -- (262.87,217) -- (386.58,232.08) ;

\draw [color = darkgray, dashed] (182.58,50) -- (182.58,87.08);
\draw [color = darkgray, dashed] (386.58,232.08) -- (454.5,232.08);

%Draw line: 
\draw[color = black] (140,96.83) -- (325,277);
\draw (320,280) node [anchor= east]    {$(\lambda^i)^\top  C$};

%Draw segment:
\draw[color = blue] (211.11,166.83) -- (262.87,217);

%Shape: Nodes
\node[vertex] at (182.58,87.08) {};
\node[vertex2] at (195.42,125.75) {};
\node[vertex] at (211.11,166.83) {};
\draw (211.11,160) node [anchor= west]    {$y^i$};
\node[vertex] at (262.87,217) {};
\draw (262.87,205) node [anchor= west]    {$y^{i+1}$};
\node[vertex] at (386.58,232.08) {};
\node[vertex2] at (250,204.5) {};
\node[vertex2] at (225,180) {};
\node[vertex2] at (360,228.8) {};

\draw (200,185) node [anchor=north west]    {$F_i$};
% Text Node
\draw (110,43.4) node [anchor=north west]    {$c^{2}(f)$};
% Text Node
\draw (433,265) node [anchor=north west]    {$c^{1}(f)$};

\draw (330,67) node [anchor=north west] {$\text{conv}(\mathcal{Y})+\R^2_{\geq}$};

\end{tikzpicture}
\caption{  Illustration of two neighboring extreme points \(y^i\) and \(y^{i+1}\), the maximally nondominated edge $F_i$ in blue and the associated weight vector $\lambda^i$.}
\label{Fig: segmentOnConvexHull}
\end{figure}

Then $f^i$ and $f^{i+1}$ are both optimal flows for the single-objective weighted-sum (MCF) program $(P_{\lambda^i})$~\cite{Eusebio09}. Hence, determining all optimal solutions for $(P_{\lambda^i})$ gives all supported efficient flows whose image lies in between $F_i$.
% We now want to determine all supported flows which images lie on the line segment of the convex hull between the extreme points $y^i $ and $ y^{i+1}$.  We define these maximally nondominated edge as $F_i$.
\Cref{Fig: segmentOnConvexHull} illustrates the objective function of the weighted-sum problem $(P_{\lambda^i})$ and the maximally nondominated face between two consecutive extreme points in the outcome space.

\begin{algorithm}[h!]
\KwData{ $(D,l,u,b,c^1,c^2)$ }
\KwResult{ The  set of all supported flows $\mathcal{X}_S$}
\(\mathcal{X}_S =  \varnothing\)\;
\tcp{Determine all $N$ extreme supported points~$y^i$ and for each one corresponding extreme flow $f^i$, sorted non-decreasingly in \(c^1(f)\)}
$\{(y^i,f^i) \colon i\in \{1,\dots,N\}\} = $ EnhancedParametricNetworkAlgortihm($D$) \;

\For{$i=1,\dots,N-1$}{
    \(\lambda_1^i = c^2(f^i)-c^2(f^{i+1}) \);\quad  \(\lambda_2^i = c^1(f^{i+1})-c^1(f^i) \) \;  
	%$D_{i,{i+1}} \leftarrow$ the network with cost function w.\,r.\,t.\ $\lambda^i$ \;
	\(\mathcal{X}_S = \mathcal{X}_S\, \cup\,\)FindAllOptimalFlows($P_{\lambda^i},f^i$) \;
	\tcp{Return only flows with $c^1(f)\neq c^1(f^i)$ to avoid repetitions.}
    
}
\Return $\X_S$.
\BlankLine
\caption{FindAllSupportedEfficientFlowsBiObjective}\label{algo:allSupportedEfficientFlows}
\end{algorithm}

\begin{thm}
Given the directed network $(D,l,u,b,(c^1,c^2)^\top)$,~\Cref{algo:allSupportedEfficientFlows} determines the set of all supported flows \(\mathcal{X}_S\) in $\mathcal{O}(Nn(m+n\log n) + \mathcal{S} (m+n))$ time.
\end{thm}
\begin{proof} 
Any supported point must lie at least on one maximally nondominated edge, and only extreme supported points $y^i$ for $i=2,\dots,N-1$  lie in two maximally nondominated edges, namely $F_{i}$ and $F_{i+1}$. According to~\Cref{thm:supp_face}, all supported flows can be determined as optimal flows for the different weighted-sum problems $(P_{\lambda^i})$. 
Moreover, unsupported flows correspond to suboptimal solutions for all weighted-sum problems. Since we only store flows for $i=2,\dots,N-1$ where $c^1(f)\neq c^1(f^i)$, no supported flow is stored twice. Thus,~\Cref{algo:allSupportedEfficientFlows} determines the complete set of all supported flows. 

The enhanced parametric network approach~\cite{raith17} requires $\mathcal{O}(Nn(m+n\log n)+M)$ time, where $M$ is the time required to solve a single-objective minimum cost flow problem. Since the algorithm determines the extreme points in a decreasing order of $c^1(f)$, we do not need additional time to sort the extreme points. Defining the weight vectors $\lambda^i$ for $i=1,\dots,N-1$ and building the network with the corresponding cost function takes $\mathcal{O}(N(n+m))$ time. Determining all $\mathcal{F}_i$ optimal flows for one weighted-sum problem $(P_{\lambda^i})$ using the algorithm FindAllOptimalFlows($P_{\lambda^i}$,$f^i$) from~\cite{KONEN2022333}, requires $\mathcal{O}(\mathcal{F}_i(m+n)+m\,n)$ time, where $f^i$ is a corresponding optimal solution to $F_i$. Since the image of every supported efficient solution lies at most on two maximally nondominated faces, it holds that $\sum_{i=1}^{N-1} \mathcal{F}_i < 2 \,\mathcal{S} $. We must consider $N-1$ of these single-objective minimum cost flow problems corresponding to the maximally nondominated faces. Hence, \Cref{algo:allSupportedEfficientFlows} requires overall $\mathcal{O}(Nn(m+n\log n) + \mathcal{S} (m+n) )$ time.
\end{proof}

Note that the determination of all supported efficient flows could easily be integrated during the enhanced parametric network approach~\cite{raith17}. Whenever a new extreme nondominated vector is found, determine all optimal flows to $(P_{\lambda^i})$ with the AOF~algorithm.

\subsection{Multi-Objective Minimum Cost Flow Problems}
\label{sec: multi}
In the following, we derive an algorithm to determine the complete set of all supported efficient flows, and hence, all supported nondominated vectors for the multi-objective integer minimum cost flow problem.

First, we need to determine the set of extreme supported nondominated vectors and the associated weight space decomposition. We can determine all extreme nondominated vectors using the dual Benson's  algorithm~\cite{Ehrgott2012}. However, since the lexicographic version of a MOIMCF problem can be solved in polynomial time~\cite{Hamacher2007}, we may also use the lexicographic dual Benson's algorithm recently presented by Bökler and Mutzel in~\cite{Boekler2015}. Both versions work with the upper image $\mathcal{P}\coloneqq {\conv(\mathcal{Y}+\R^d_{\geq})}$ and its dual polyhedron, or \emph{lower image} $\mathcal{D}$. 

While we work with normalized weight vectors $\lambda\in \Lambda_d$, it suffices to consider the so-called \emph{projected weight space} $\{(\lambda_1,\ldots, \lambda_{d-1})\in \Lambda_{d-1}\colon \sum_{i=1}^{d-1} \lambda_i<1 \}$ and calculate the normalized weighting vector 
$\ell(v)\coloneqq (v_1,\dots, v_{d-1},\sum_{i=1}^{d-1} v_i)$ of a projected weight \(v\)
%$\lambda_d = 1- \sum_{i=1}^{d-1} \lambda_i$ 
when needed. %Let $\ell(v)\coloneqq (v_1,\dots, v_{d-1},\sum_{i=1}^{d-1} v_i)$ . 
The dual problem \eqref{eq:DPl} of the weighted sum-problem \eqref{eq:P_l} is given by 
\begin{equation}\label{eq:DPl}\tag{$D_{\lambda}$}
    \begin{array}{rl}
      \max & b^\top u \\
      \text{s.t.} & A^\top u= C^\top \lambda\\
      &u \in \R_{\geq}^{m}.
    \end{array}
\end{equation}
The dual polyhedron $\mathcal{D}$ then consists of vectors $(\lambda_1, \dots, \lambda_{d-1}, b^\top u)$ with  $\lambda\in \Lambda_d^0$ and solutions $u$ of \eqref{eq:DPl}. 
Following the duality theory of polyhedra, there exists a bijective mapping $\Psi$ between the set of all faces of $\mathcal{P}$ and the set of all faces of $\mathcal{D}$ such that $\Psi$ is \emph{order reversing}, i.\,e., if two faces $F_1$ and $F_2$ of
$\mathcal{P}$ satisfy $F_1\subseteq F_2$ then $\Psi(F_1) \supseteq \Psi(F_2)$ and $\Psi(F_1)$ and $\Psi(F_2)$ are faces of $\mathcal{D}$, see e.g., \cite{Schulze2019}. Thus, an extreme point of $\mathcal{D}$ corresponds to a facet of $\mathcal{P}$, and an extreme point of $\mathcal{P}$ corresponds to a facet of $\mathcal{D}$. The dual Benson's algorithm solves a MOLP by computing the extreme points of $\mathcal{D}$. For more details on the dual Benson's algorithm or its lexicographic version, we refer to~\cite{Boekler2015,Ehrgott2012}.

Thus, we obtain all extreme nondominated vectors and one corresponding extreme efficient solution for each of the extreme nondominated vectors, as well as all facets  of $\mathcal{P}$. On this basis, we yield the \emph{weight space decomposition} using the dual (lexicographic) Benson's algorithm.   
The set of weighting vectors associated with a point $y\in \mathcal{Y}$ is given by 
\[
\mathcal{W}(y)\coloneqq \left\{w \in \Lambda^0_{d}\colon w^\top  y \leq w^\top  y^{\prime} \text{ for all } y^{\prime} \in \mathcal{P}\right\}. 
\]

Note that the facets of $\mathcal{P}$ may only be \emph{weakly} nondominated, i.\,e., they might contain dominated (integer feasible) vectors. Recall that all supported nondominated vectors can be determined by a parametric MCF problem \eqref{eq:P_l} %programs in the form of ${P}(\lambda)\coloneqq  \min \{\lambda^\top  c(f)\colon f \in \mathcal{X}\}$ 
for some weight vector $\lambda\in \Lambda_d$. However, the weight vectors corresponding to the facets of $\mathcal{P}$ might have components equal to zero $\lambda_i=0$ for $i\in\{1,\dots,d\}$, i.\,e., $\lambda \in \Lambda_d^0$.
In the following, we describe a recursive algorithm to obtain the weight vectors for all maximally nondominated faces.

Let $U$ be the set of all extreme points in the lower image $\mathcal{D}$ and $\{\lambda^u\colon u\in U\}$ the set of corresponding weight vectors. Let, furthermore, 
be $F_u$ the facet of $\mathcal{P}$ corresponding to \(u\in U\). Then, we call two extreme points \(u\) and \(u'\) of \(\mathcal{D}\) \emph{adjacent}, iff 
$\dim(F_u \cap F_{u'})=d-2$. In the following, we will denote the set of adjacent extreme points for $u\in U$ by \(Q_u\subseteq U\).

Recall that the intersection of 
$k$ adjacent facets yield a $d-k$ dimensional face. For each $\lambda^u \in \Lambda_d$ (i.\,e., $\lambda^u>0$) we know that all vectors on the facet $F_u$ are supported nondominated vectors. Thus, we only have to solve the all-optimum flow problem on $(P_{\lambda^u})$. 
Since some solutions may lie in the same sub-faces of adjacent facets, we have to ensure that no solution is stored twice. In order to do so, we keep track of the neighbouring extreme points during Benson's algorithm and store all already processed adjacent extreme points of $u\in U$ in a list $\delta_u$. 

There may exist maximally nondominated faces (with dimension less than \(d-1\)), which are intersections of a number of facets for which the corresponding weight vector equals zero in at least one component. We call these facets \emph{weakly nondominated facets} of $\mathcal{P}$. In order to determine all supported efficient solutions, we investigate  nondominated faces which are intersections of weakly supported facets. %Note that the upper image $\mathcal{P}$ does not contain any redundant facets. 
With $U_>\coloneqq \{u\in U \colon \lambda^u > 0 \}$ we denote the set of extreme points of $\mathcal{D}$ corresponding to nondominated facets and with 
%\(U_0\coloneqq  U \backslash U_>\) the set of extreme points corresponding to weakly nondominated but not nondominated facets.}
\(U_0\coloneqq \{u\in U \colon \lambda^u \geq 0 \}\) the set of extreme points corresponding to all weakly nondominated facets (\(U_>\subseteq U_0\)).
Note that weakly nondominated faces can contain supported nondominated points only at its (relative) boundary, while unsupported nondominated points can be located also in its (relative) interior.

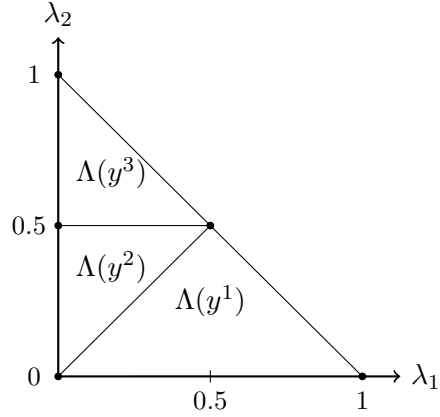
\begin{figure}[t]
\centering
\begin{tikzpicture}
\tikzstyle{vertex}=[circle,fill=black,draw=black,minimum size=2.5pt,inner sep=0]

\draw [<->,thick] (0,4.5) node (yaxis) [above] {$\lambda_2$}
        |- (4.5,0) node (xaxis) [right] {$\lambda_1$};
        
\draw (0,4) node (a_1) [xshift=-9pt] {\small$1$} -- (4,0) node (a_2) [yshift=-9pt] {\small$1$};
\draw (0,2) node (a_2) [xshift=-11pt] {\small$0.5$} -- (2,2) coordinate (a_3);
\draw (0,0) node (a_4) [xshift=-9pt] {\small$0$} -- (2,2) coordinate (a_5);
\draw (2,0) node (a_6) [yshift=-9pt] {\small$0.5$} -- (1,0) coordinate (a_7);
\draw (2,-0.08)--(2,0.08);

\node[vertex] at (0,0) (int1) {};
\node[vertex] at (4,0) (int1) {};
\node[vertex] at (0,4) (int1) {};
\node[vertex] at (2,2) (int1) {};
\node[vertex] at (0,2) (int1) {};

\node at (0.7,2.7) (int1) {$\Lambda(y^3)$};
\node at (0.7,1.45) (int1) {$\Lambda(y^2)$};
\node at (2,1) (int1) {$\Lambda(y^1)$};
\end{tikzpicture}
 \caption{Weight space decomposition to the upper image of~\Cref{obj-ex1}. Here it holds $\lambda_3=1-\lambda_2-\lambda_1.$}
 \label{fig:wsd}
\end{figure}

\Cref{fig:wsd} presents the weight space decomposition for the example given in \Cref{2} and illustrated in \Cref{obj-ex1}. Any point in the weight space decomposition corresponds to an $\lambda^u$ for each $u\in U$. Here $U_>= \varnothing$. However, there do exist weights in the lines in the interior connecting two adjacent extreme points of \(\mathcal{D}\) which correspond to the maximally nondominated faces \([y^1,y^2]\) and \([y^2,y^3]\).

\begin{algorithm}[h!]
\KwData{ $(D,l,u,b,c)$ }
\KwResult{ The complete set of all \strictly supported efficient flows }
\BlankLine

$\{U,Q_u, F_u ,\lambda_u, f^u \colon u\in U\} =$ BensonLex$(D,l,u,b,c)$ \;
\tcp{Determine all extreme points of $\mathcal{D}$, the corresponding facets of $\mathcal{P}$, and weight vectors.}
$\delta_u=\varnothing \quad \forall u\in U$\;
\For{$u\in U_>$}{
    $ \X_S = \X_S \, \cup $ FindAllOptimalFlows$(P_{\lambda^u},f^u)$\; 
    \tcp{In the FindAllOptimalFlows algorithm only store flows $f$ for which 
    $\langle\lambda^u,C\,f\rangle
    \neq  \min\{\langle\lambda^{u'},C\,f'\rangle \colon f' \in \mathcal{X}\} \text{ for any } u'\in \delta_u $}
    \For{$u'\in Q_u$}{$\delta_{u'}= \delta_{u'} \cup \{u\}\;$}
}
$w_u = \{ {\lambda^{u'}} \colon u'\in \delta_u \}  \quad \forall u\in U$\;
$B = \varnothing$\;

\For{ $u \in U_0\backslash U_>$}{
    $\Tilde{U} = \{ u \}$ \; 
     $ \X_S = \X_S \, \cup $ ConsiderSubFaces$(\Tilde{U},U_0\backslash U_>, B, Q_u,\delta_u,\lambda^u,f^u,w_u \; \forall u \in U)$\;
    $B = B \cup \{u\}$
}
\Return $\X_S$
\BlankLine
\caption{FindAllSupportedEfficientFlows}\label{algo:allSupportedEfficientFlowsMulti}
\end{algorithm}

\begin{algorithm}[h!]
\KwData{ $(\Tilde{U},U_0\backslash U_>,B,Q_u,\delta_u,\lambda^u,f^u,w_u) \quad \forall u\in U$ }
\KwResult{  The set $\bar{\X}$ of all \strictly supported efficient flows in the maximally nondominated sub-faces of the facet $F_u$ of $\mathcal{P}$, which are not lying in an already investigated faces}
\BlankLine

Let $\bar{U}\coloneqq \{  u^{\prime} \in (\bigcap_{u\in \Tilde{U}} Q_u\cap U_0\backslash U_>)\backslash (B \cup \Tilde{U}) \}$\;
$\bar{\X} = \emptyset$ \;

\For{ $u\in \bar{U}$}{
$\Tilde{U}= \Tilde{U}\cup \{u\}$\;
    \If{$|\Tilde{U}| \le d-1 $}
    {
    $\lambda^{k}\coloneqq  \sum_{i=1}^{|\Tilde{U}|} \ell_i \,\lambda^{\Tilde{u}_i}$ for an $\ell \in \{ \sum_{i=1}^{|\Tilde{U}|} \ell_i =1, \ell_i > 0   \quad \forall i \in \{1,\ldots, |\Tilde{U}|\}\}, \Tilde{u}_i\in \Tilde{U} $\;
    
    \If{$\lambda^k > 0$}
        {
        $\bar{\X} = \bar{\X}  \, \cup $ FindAllOptimalFlows$(P_{\lambda^k},f^{\tilde{u}_1})$\; 
        \tcp{Only store flows $f$ for which 
        $\langle\lambda^k,Cf\rangle\neq \{ \min\{ \langle\lambda^{u'},Cf'\rangle \colon f' \in \mathcal{X}\} \colon \lambda^{u^{\prime}} \in \bigcup_{u\in \Tilde{U}} w_u \}$}
        \For{$u' \in \cup_{u\in \Tilde{U}} Q_u $}
            {
                $w_{u'}= w_{u'} \cup \{\lambda^k\}$\;
            }   
        }
    \Else{ $\bar{\X} = \bar{\X}  \, \cup $ ConsiderSubFaces$(\Tilde{U},U_0\backslash U_>, B, Q_u,\delta_u,\lambda^u,f^u,w_u \; \forall u \in U)$\;}
}
$B = B\cup \{u\}$
}
\Return $\bar{\X}$
\BlankLine
\caption{ConsiderSubFaces}\label{algo:considerSubFaces}
\end{algorithm}

\begin{thm}\label{thm:outputpoly}
Given the directed network $(D,l,u,b,c)$,~\Cref{algo:allSupportedEfficientFlowsMulti} determines the complete set of all \strictly supported flows in $\mathcal{O}(N^{\lfloor \frac{d}{2}\rfloor}( \operatorname{poly}(n,m) + N\log N + N^{\lfloor \frac{d}{2}\rfloor}\mathcal{S}( m+n + N^{\lfloor \frac{d}{2}\rfloor}d))$ time.  
\end{thm}
\begin{proof}
\textit{Correctness:} Any \strictly supported efficient flow must lie in at least one maximally efficient face. However, a \strictly supported efficient flow can lie in multiple face. Due to~\Cref{thm:supp_face}, all \strictly supported flows are found by determining all optimal flows for each weight vector $\lambda^i$ corresponding to a maximally nondominated face. 
Moreover, no weakly supported flow can be optimal for a parametric network flow problem with one of these cost functions. While the algorithm iterates through all maximally nondominated faces, only flows are stored that have not been considered yet. Thus,~\Cref{algo:allSupportedEfficientFlowsMulti} determines the complete set of all supported efficient solutions. 

\textit{Run-time:} Benson's Algorithm requires $\mathcal{O}(N^{\lfloor \frac{d}{2}\rfloor}(\operatorname{poly}(n,m)+ N \log N))$ time~\cite{Boekler2015}. Thereafter, we consider each face at most once. The number of all faces can be bounded by $\mathcal{O}(N^{\lfloor \frac{d}{2}\rfloor})$~\cite{Boekler2015}.  We check for each weakly nondominated face if a strict convex combination with adjacent weight vectors yields a weight vector $\lambda>0$ componentwise strictly greater than zero. In this case, we call~\Cref{algo:considerSubFaces}. Note that $\lambda$ is not part of the input of \Cref{algo:allSupportedEfficientFlowsMulti}. However, it can be shown that these encoding lengths can be bounded by $\mathcal{O}(\operatorname{poly}(n,m))$~\cite{Boekler2015}. The convex combination can be computed in $\mathcal{O}(N^{\lfloor \frac{d}{2}\rfloor}d)$. First, we must create the weight vector through a strictly convex combination for each maximally nondominated face, which is not a facet. Afterwards, we solve the AOF problem for \eqref{eq:P_l} for all of these maximally nondominated faces in time $\mathcal{O}(F_i(m+n+ N^{\lfloor \frac{d}{2}\rfloor}d))+m\,n)$, where $F_i$ is the number of optimal solutions for the current weighted sum problem \eqref{eq:P_l}. 
Additionally, it takes $\mathcal{O}(N^{\lfloor \frac{d}{2}\rfloor}d)$ time to check if the flow is also optimal for an adjacent already considered maximally nondominated face. 
Since each flow may be contained in all faces, we obtain the bound $\sum F_i \leq \mathcal{O}(N^{\lfloor \frac{d}{2}\rfloor}\mathcal{S})$, where $\mathcal{S}$ is the number of all supported efficient flows.
Hence, \Cref{algo:allSupportedEfficientFlowsMulti} requires overall $\mathcal{O}(N^{\lfloor \frac{d}{2}\rfloor}( \operatorname{poly}(n,m) + N\log N + N^{\lfloor \frac{d}{2}\rfloor}\mathcal{S}( m+n + N^{\lfloor \frac{d}{2}\rfloor}d))$ time.  
\end{proof}

\section{Conclusion}\label{chapt:concl}

This paper concludes that there is no output-polynomial algorithm for a MOIMCF problem with a fixed number of $d$ objectives that determines all weakly supported nondominated vectors unless $\mathbf{P} = \mathbf{N P}$. In contrast, this paper presents output-polynomial time algorithms for determining all supported efficient solutions for BOIMCF problems and general MOIMCF problems with a fixed number of objectives.
First, the approach determines all extreme supported nondominated vectors and the weighting vectors for each maximally nondominated face. Then, it successively determines all supported efficient solutions in the preimage of each maximally nondominated face by determining all optimal solutions for the corresponding single-objective parametric network flow problem using the all optimum flow algorithm recently presented in~\cite{KONEN2022333}. 

However, it might be that many supported efficient flows may be mapped to the same vector in the objective space.  

Thus, often, a minimal complete set (all nondominated vectors and one (efficient) preimage for each of them) is considered as a solution to a multi-objective optimization problem~\cite{serafini87some}. An open question remains if  an output-polynomial time algorithm exists to determine all supported nondominated vectors for a MOIMCF problem with a fixed number of objectives.  

Even though an output-polynomial time algorithm to determine all nondominated vectors for MOIMCF problems does not exist, even for the bi-objective case~\cite{Boekler2017}, future research could focus on new approaches to compute also unsupported nondominated vectors/efficient solutions in bi- or even multi-objective MCF problems. Unsupported solutions may be good compromise solutions and should thus not be neglected completely. Note that the difficulty to compute unsupported solutions is not a specific property of multi-objective integer network flow problems but arises in many integer and combinatorial optimization problems and is one reason for their computational complexity, in general~\cite{ehrgott00hard, figueira17easy}. One way to overcome this computational burden---at least to a certain degree---could be to determine unsupported solutions only in regions of the Pareto front that are not well represented by the set of supported nondominated points.

\bigskip

\textbf{Acknowledgements.} We would like to thank Fritz Bökler for the very kind and in-depth discussion on the properties of the dual Benson's algorithm and on the complexity of multi-objective combinatorial optimization problems. We thank Kathrin Klamroth for many fruitful discussions on our work. David Könen acknowledges financial support from Deutsche Forschungsgemeinschaft within the
project number 441310140.

\end{document}